\documentclass{article}
\usepackage{theorem,amssymb,amsfonts}
\usepackage{variations}
\date{}
\newtheorem{lemma}{Lemma}[section]   \newtheorem{corollary}{Corollary}[section]
\newtheorem{theorem}{Theorem}[section]

\def\subjclass{{\bf 2000 Mathematics Subject Classification:}}
\def\keywords{{\bf Key Words and Phrases:}} 
\begin{document}
\title{\bf Study of some semi-linear elliptic equation}
\author{\bf Anouar Ben Mabrouk\footnote{Department of Mathematics, Faculty of Sciences, Monastir, Tunisia. Email: anouar.benmabrouk@issatso.rnu.tn}}
\maketitle
\begin{abstract}
We propose in this paper to study nodal solutions of some nonlinear elliptic equations derived from the famous equation of Brezis-Nirenberg and we analyse in some cases the possible singularities of radial solutions at the origin.\\
{\keywords} elliptic equations, nodal solutions, Sturm comparison theorem, Brezis-Nirenberg equation.\\ {\subjclass} 35J25, 35J60.
\end{abstract}
\section{Main Results}
In this paper we focus on the study of the equation
\begin{equation}\label{mainproblem}
\Delta\,u+u-|u|^{-2\theta}u=0,\,\hbox{in}\,\mathbb{R}^d,
\end{equation}
with $d>1$, and $0<\theta<\frac{1}{2}$. We focus essentially on radial solutions. The radial version of problem (\ref{mainproblem}) provided with the value of the solution $u$ at the origin is
\begin{equation}\label{mainproblemradial}
\left\{\matrix{u''+\frac{d-1}{r}u'+u-|u|^{-2\theta}u=0&,&r\in(0,+\infty),\hfill\cr
\medskip\,u(0)=a&,&u'(0)=0,\hfill}\right.
\end{equation}
where $a\in\mathbb{R}$. In the rest of the whole paper we denote 
$$
g(s)=1-|s|^{-2\theta},\;\;f(s)=sg(s)=s-s|s|^{-2\theta}\,\hbox{and}\,F(s)=\displaystyle\frac{s^2}{2}(1-\displaystyle\frac{|s|^{-2\theta}}{1-\theta}).
$$ 
We propose to recall some results already known or easy to handle on the study of problem (\ref{mainproblem}) or one of the problems related and cited above. To do this we recall the properties of $g$, $f$ and $F$. Because of the parity properties of these functions, we only provide their variations on $(0,+\infty)$.
$$
\begin{variations}
s&0&&1&&&\pI\\
\hline
g'(s)&\bb&+&&\ga+\\
\hline
\m {g(s)}&\mI&\c&\h\z&\z&\c&\h1\\
\hline
{g(s)}&0&\ga-&0&\ga+\\
\end{variations}
$$
$$
\begin{variations}
s&0&&s_{\theta}&&1&&\ga\pI\\
\hline
f'(s)&\bb&-&\z&+&&\ga+&\\
\hline
\m {f(s)}&\h\z&\d&f(s_{\theta})&\c&\h\z&\z&\c&\h\pI\\
\hline
{f(s)}&0&\ga-&\ga0&\ga+&\\
\end{variations}
$$
$$
\begin{variations}
s&0&&1&&p&&\ga\pI\\
\hline
F'(s)&\z&-&\z&+&&\ga+&\\
\hline
\m {F(s)}&\h\z&\d&F(1)&\c&\h\z&\z&\c&\h\pI\\
\hline
{F(s)}&0&\ga-&\ga0&\ga+&\\
\end{variations}
$$
The parameter $p=\displaystyle\frac{1}{(1-\theta)^{\frac{1}{\theta}}}$ is the unique real number in $(1,+\infty)$ such that $F(p)=0$. We recall finally that we shall use many times the energy of the solution $u$ defined for $r\geq0$ by
$$
E(r)=\displaystyle\frac{1}{2}u'^2(r)+\displaystyle\int_0^{u(r)}sg(s)ds=\displaystyle\frac{1}{2}u'^2(r)+\displaystyle\int_0^{u(r)}f(s)ds=\displaystyle\frac{1}{2}u'^2(r)+F(u(r)).
$$
The first result is stated as follows.
\begin{theorem}\label{theorem1}
The solution $u$ of problem (\ref{mainproblemradial}) is oscillating around $1$ or $-1$ for any $a\in]-1,1[\setminus\{0\}$ with no zeros in $(0,\infty)$.
\end{theorem}
Next we study the case where the origin value $u(0)=a$ is not in the $\pm1$-attractive zone. We prove that there are also different zones to be distinguished. We obtained the following result.
\begin{theorem}\label{theorem2}
\begin{enumerate}
\item[i.] For $1<a<p$, the solution $u$ of problem (\ref{mainproblemradial}) ......
\item[ii.] For $a>p$, the solution $u$ of problem (\ref{mainproblemradial}) ......
\end{enumerate}
\end{theorem}
The following result deals with the existence and uniqueness of the solution.
\begin{corollary}\label{corollarysolutionoscillenteaentre1etp2}
For all $a\in]1,p[$, problem (\ref{mainproblemradial}) has a unique solution $u>0$ which is oscillatory around $1$.
\end{corollary}
\section{On the existence and uniqueness of solutions}
\begin{lemma}\label{lemmeu1alorssaderiveenot0}
For all $a\in(0,p)$, the solution $u$ of (\ref{mainproblemradial}) satisfies the assertion
$$
u(\zeta)=0,\,\hbox{\it for some}\,\zeta\Longrightarrow\,u'(\zeta)\not=0,
$$
except if $u\sim0$.
\end{lemma}
\begin{lemma}\label{lemmesolutionpositive1p2}
For all $a\in]0,p[$, with $p=\displaystyle\frac{1}{(1-\theta)^{\frac{1}{2\theta}}}$, problem (\ref{mainproblemradial}) has a unique positive solution $u$.
\end{lemma}
Indeed, denote for $r\in(0,+\infty)$ and consider the system
\begin{equation}\label{systemiteratif}
\left\{\matrix{u(r)=a+\displaystyle\int_0^rv(s)ds,\hfill\cr
\medskip\,v(r)=-\displaystyle\frac{1}{r^{d-1}}\displaystyle\int_0^rs^{d-1}u(s)g(u(s))ds.\hfill}\right.
\end{equation}
Using standard arguments from iterative methods in functional analysis, we observe that such a system has a unique local solution $(u,v)$ on $r\in(0,\delta)$ for $\delta>0$ small enough. The solution satisfies $u(0)=a$, $v(0)=0$. Furthermore, $u>0$, $v<0$, and $u$ and $v$ are $\mathcal{C}^2$ on $]0,\delta[$ and
$$
u'(r)=v(r)\quad\hbox{and}\quad\,v'(r)=-\displaystyle\frac{d-1}{r}v(r)-u(r)g(u(r)).
$$
We now study the differentiability at $0$. Using L'Hospital rule, we obtain
$$
u''(0)=v'(0)=\displaystyle\lim_{r\rightarrow0}\displaystyle\frac{v(r)}{r}=-\displaystyle\frac{ag(a)}{d}.
$$
On the other hand,
$$
\matrix{\displaystyle\lim_{r\rightarrow0}v'(r)&=&\displaystyle\lim_{r\rightarrow0}u''(r)\hfill\cr
\medskip&=&-\displaystyle\lim_{r\rightarrow0}\hfill\cr
\medskip&=&\left[(d-1)\displaystyle\frac{v(r)}{r}+u(r)g(u(r))\right]\hfill\cr
\medskip&=&\displaystyle\left[(d-1)\displaystyle\frac{ag(a)}{d}+ag(u(a))\right]\hfill\cr
\medskip&=&-\displaystyle\frac{ag(a)}{d}.\hfill}
$$
Hence, $u$ is $\mathcal{C}^2$ at $0$. It suffices then to prove that $u>0$ on $(0,+\infty)$ to guarantee the existence and uniqueness on $(0,+\infty)$. We suppose by contrast that $u(\zeta)=0$ for some $\zeta>0$. The evaluation of the energy $E$ gives
$$
E(\zeta)=\frac{1}{2}u'^2(\zeta)<E(0)=F(a)<0
$$
because of the fact $1<a<p$. Which leads to a contradiction.\\

Let $u$ be a compactly supported solution of problem (\ref{mainproblemradial}) already with $a>1$ and let $R=\inf\{r\in(0,\infty),\;u(s)=0,\;\forall\,s\geq\,r\}$. Henceforth, $u$ is a solution of the problem
\begin{equation}\label{mainproblemB0R}
\left\{\matrix{\Delta\,u+u-|u|^{-2\theta}u=0&in&B(0,R),\hfill\cr
\medskip\,u=0&on&\partial\,B(0,R).\hfill}\right.
\end{equation}
Recall that it is well known that for $R<\sqrt{\lambda_1(B(0,1))}$ the first eigenvalue of $-\Delta$ on the unit ball, problem (\ref{mainproblemB0R}) has no positive solution. See \cite{chen1}, \cite{davilamontenegro} and the references therein. Consequently, we will assume for the rest of this part that $R\geq\sqrt{\lambda_1(B(0,1))}$ and consider the radial expression of (\ref{mainproblemB0R}),
\begin{equation}\label{mainproblemradialB0R}
\left\{\matrix{u''+\frac{d-1}{r}u'+u-|u|^{-2\theta}u=0&,&r\in(0,+\infty),\hfill\cr
\medskip\,u(R)=0.\hfill}\right.
\end{equation}
We will discuss the behavior of the solution $u$ relatively to the values $u'(R)$. Two situations can occur. First, $u'(R)<0$. It results that $u'(r)<0$ on a small interval $(R-\varepsilon,R+\varepsilon)$. Therefore, $u(r)<0$ on $(R,R+\varepsilon)$ which contradicts the definition of $R$. Next, for $u'(R)=0$, we get $E(R)=0$

Denote for the rest of the paper $\rho_a$ the first zero of the solution $u$ of problem (\ref{mainproblemradial}) for $a>p$. We have
\begin{lemma}\label{lemmeasuperieurap2}
For all $a>p$, $\rho_a<\infty$.
\end{lemma}
\ {\it Proof.} Suppose $u$ a solution of problem (\ref{mainproblemradial}) with $a>p$ and $\rho_a=\infty$. The solution $u$ starts as decreasing from $a=u(0)$. Suppose that it remains decreasing on its whole domain $(0,+\infty)$. Thus it has a limit $L$ as $r\rightarrow+\infty$. Thus $L=0$ or $L=1$. For $L=0$ and $r$ large enough, we obtain $u(r)=A\cos(r)+B\sin(r)$ which is contradictory. The case where $L=1$ is analogous. Consequently $\rho_a<+\infty$.\\

We now study the behavior of the solution $u$ on the whole domain $(0,+\infty)$. Denote $r_0$ the first critical point of the solution $u$ of problem (\ref{mainproblemradial}) with $a>p$. There are four possible situations. The case $u(r_0)>1$ with equation (\ref{mainproblemradial}) implies that
$$
0=\displaystyle\int_0^{r_0}\bigl(s^{d-1}u'(s)\bigr)'ds=-\displaystyle\int_0^{r_0}s^{d-1}u(s)g(u(s))ds<0
$$
which is impossible. Next, for $u(r_0)=1$, the solution $u$ will be a solution of the problem
\begin{equation}\label{mainproblemradialur0egal1}
\left\{\matrix{u''+\frac{d-1}{r}u'+u-|u|^{-2\theta}u=0&,&r\in(r_0,+\infty),\hfill\cr
\medskip\,u(r_0)=1&,&u'(r_0)=0.\hfill}\right.
\end{equation}
Therefore, $u\equiv1$, for any $r\geq\,r_0$, which is contradictory by the same argument as above. We now assume that $0<u(r_0)<1$. Implying Theorem 1.2 in \cite{anouarmohamednodal}, we observe that $u$ is oscillating around $0$ with no zeros. Which is contradictory. Now we examine the last case $u(r_0)=0$. In this case, we obtain $u(r_0)=u'(r_0)=0$ and $\rho_a>\sqrt{\lambda_1(B(0,1))}$.

\section{Proof of Theorem \ref{theorem1}.}
We recall firstly that some situations which are somehow more general are developed in \cite{anouarmohamednodal}. The proof developed here is inspired from there. Let $a\in(0,1)$ and $u$ be the solution of problem (\ref{mainproblemradial}). It holds that $u''(r)>0$ on a small interval $(0,\varepsilon)$ for $\varepsilon$ small enough positive. Consequently, $u'$ is strictly increasing on $(0,\varepsilon)$. Which yields that $u'(r)>0$ on $(0,\varepsilon)$. Thus $u$ is strictly increasing on $(0,\varepsilon)$ for $\varepsilon$ small enough positive. So that, $u(r)>a$ on $(0,\varepsilon)$. We will prove that the value $a$ is taken only for $r=0$. Indeed, suppose not, and let $\zeta>0$ be the first point satisfying $u(\zeta)=a$. The evaluation of the energy $E(r)$ at $0$ and $\zeta$ yields that
$$
E(0)=F(a)>E(\zeta)=\displaystyle\frac{1}{2}u'^2(\zeta)+F(a)
$$
which is contradictory. So, the solution $u$ starts increasing with origin point $u(0)=a$ and did not reach it otherwise. We next prove that it can not continue to increase on its whole domain $(0,+\infty)$. Suppose contrarily that it is increasing on $(0,+\infty)$ and denote $L$ its limit as $r\rightarrow+\infty$. Of course, such a limit can not be infinite because of the energy of the solution. Next, the finite limit is a zero of the function $f(s)$. Therefore, $L=1$. But, this yields $u''(r)>0$ as $r\rightarrow+\infty$ (Recall that $f(s)<0$ on $(0,1)$). In the other hand, equation (\ref{mainproblemradial}) guaranties that $\frac{u''(r)}{u'(r)}\sim-\frac{d-1}{r}<0$ as $r\rightarrow+\infty$ which means that $u''(r)<0$ as $r\rightarrow+\infty$ leading to a contradiction. We therefore conclude that $u$ is oscillatory. Let $t_1$ be the first point in $(0,+\infty)$ such that $u'(t_1)=0$. It holds that $u(t_1)>1$. If not, by multiplying equation (\ref{mainproblemradial}) by $r^{d-1}$ and integrating from $0$ to $t_1$ we obtain $0=-\displaystyle\int_0^{t_1}r^{d-1}f(u(r))>0$ which is contradictory. Thus, $u$ crosses the line $y=1$ once in $(0,t_1)$ leading to a unique point $r_1\in(0,t_1)$ such that $u(r_1)=1$. Next, using similar techniques, we prove that $u$ can not remain greater than $1$ in the rest of its domain. (Consider the same equation on $(t_1,+\infty)$ with initial data $u(t_1)$ and $u'(t_1)$). Consequently we prove that there exists unique sequences $(t_k)_k$ and $(r_k)_k$ such that
\begin{equation}\label{uzetakrk}
r_k<t_k<r_{k+1},\quad\,u(r_k)=1,\quad\,u'(\zeta_k)=0,\;k\geq1.
\end{equation}
Next, observing that $E$ is decreasing as a function of $r$, we deduce that the sequence of maxima $(u(t_k))_k$ goes to $1$ and therefore $u$.
\section{Proof of Theorem \ref{theorem2}.}
The proof is based on a series of preliminary results. We recall first that it suffices to study the case $a>0$ due to the parity properties of the function $g$ and/or $f$.
\begin{lemma}\label{lemma2.1}
For $a>1$, the solution $u$ satisfies $\bigl(u(r)<a,\;\forall\,r>0\bigr)$.
\end{lemma}
\ {\it Proof.} From equation (\ref{mainproblemradial}), we obtain $du''(0)=-ag(a)<0$. Consequently, $u''(r)<0$ for $r\in(0,\varepsilon)$ for some $\varepsilon>0$ small enough. Thus, $u'$ is decreasing strictly on $(0,\varepsilon)$ and then, $u'(r)<0$ for $r\in(0,\varepsilon)$. Therefore, $u$ is decreasing strictly on $(0,\varepsilon)$ and then, $u(r)<a$ for $r\in(0,\varepsilon)$. Let next $\zeta>0$ be the first point such that $u(\zeta)=a$, if possible. Using the energy $E$ we obtain $E(\zeta)<E(r)<E(0)$, for all $r\in(0,\zeta)$. Whenever $u'(\zeta)=0$, we obtain $E(0)<E(0)$ which is impossible. So $u'(\zeta)\not=0$, which implies that $\displaystyle\frac{1}{2}u'^2(\zeta)+E(0)<E(0)$ which is also impossible. As a conclusion, there is no positive points for which the solution $u$ reaches $a$ again.
\begin{lemma}\label{lemma2.2}
For $a>1$, the solution $u$ is not strictly decreasing on $(0,+\infty)$.
\end{lemma}
\ {\it Proof.} Assume contrarily that $u$ is strictly decreasing on $(0,+\infty)$. Thus, it has a limit $L$ as $r\rightarrow+\infty$. Two cases are possible, $L=0$ or $L=1$. We will examine them one by one.\\
\ {\it case 1.} $L=0$. Consider the dynamical system in the phase plane defined for $r\in(0,+\infty)$ by
\begin{equation}\label{mainproblemradialplandephase}
\left\{\matrix{v=u',\hfill\cr
\medskip\,v'=-\frac{d-1}{r}v+u-|u|^{-2\theta}u,\hfill\cr
\medskip\,u(0)=a&,&v(0)=0\hfill}\right.
\end{equation}
A careful study for $r\rightarrow+\infty$, yields the estimation $u\sim\,A\cos(r)+B\sin(r)$ for $r$ large enough, which is contradictory.\\
\ {\it case 2.} $L=1$. Using equation (\ref{mainproblemradial}) or (\ref{mainproblemradialplandephase}), we obtain for $r$ large enough,
$2^{d-1}v(2r)-v(r)\sim\displaystyle\frac{2^d-1}{d}r$ which leads to a contradiction.
\begin{lemma}\label{lemma2.3}
Let for $a>1$, $r_1(a)$ be the first critical point of the solution $u$ of (\ref{mainproblemradial}) in $(0,+\infty)$. Then $u(r_1)<1$.
\end{lemma}
\ {\it Proof.} Suppose not, i.e, $u(r_1)=1$ or $u(r_1)>1$. When $u(r_1)>1$, we obtain $u''(r_1)=-u(r_1)g(u(r_1))<0$. Hence, $u'$ is decreasing on $(r_1-\varepsilon,r_1+\varepsilon)$ for some $\varepsilon$ small enough. Thus, $u$ is increasing near $r_1$ at the left and decreasing near $r_1$ at the right, which is contradictory. When $u(r_1)=1$, then $u$ is a solution of the problem $u''(r)+\frac{d-1}{r}u'(r)+ug(u)=0$ on $(r_1,+\infty)$ with the initial condition $u(r_1)=1$ and $u'(r_1)=0$. Consequently, $u\equiv1$ which is contradictory.
\begin{lemma}\label{lemme2.4}
Let $a>1$ and $r_1(a)$ be the first critical point of the solution $u$ in $(0,+\infty)$. Then
\begin{description}
  \item[a.] for $r_1(a)\in]0,1[$, the solution $u$ of (\ref{mainproblemradial}) oscillates around 1, with limit 1, and thus has a finite number of zeros.
  \item[b.] for $r_1(a)\in]-1,0[$, the solution $u$ of (\ref{mainproblemradial}) oscillates around -1, with limit -1, and thus has a finite number of zeros.
\end{description}
\end{lemma}
\ {\it Proof.} In the situation {\bf a.} $u$ is a solution of the problem
\begin{equation}\label{mainproblemradialenr1(a)}
\left\{\matrix{u''+\frac{d-1}{r}u'+u-|u|^{-2\theta}u=0&,&r\in(r_1(a),+\infty),\hfill\cr
\medskip\,u(r_1(a))\in]0,1[&,&u'(r_1(a))=0.\hfill}\right.
\end{equation}
Hence, by applying Theorem \ref{theorem1}, the solution oscillates around 1, with limit 1 and thus, it has a finite number of zeros. In the situation {\bf b.} $u$ is a solution of the problem
\begin{equation}\label{mainproblemradialenr1(a)}
\left\{\matrix{u''+\frac{d-1}{r}u'+u-|u|^{-2\theta}u=0&,&r\in(r_1(a),+\infty),\hfill\cr
\medskip\,u(r_1(a))\in]-1,0[&,&u'(r_1(a))=0.\hfill}\right.
\end{equation}
Hence, for the same reasons, it oscillates around -1, with limit $-1$ and thus with a finite number of zeros.
\begin{lemma}\label{lemme2.5}
Let $a>1$ and $u$ the solution of (\ref{mainproblemradial}) in $(0,+\infty)$. The following situation can not occur.
There exists sequences $(r_k)$, $(t_k)$, $(z_k)$ and $(\zeta_k)$ satisfying
\begin{description}
\item[i.] $t_{2k-1}<z_{2k-1}<\zeta_{2k-1}<z_{2k}<t_{2k}<r_{2k}<\zeta_{2k}<r_{2k+1},\quad\forall\,k.$
\item[ii.] $u(r_k)=-u(z_k)=1,\quad\,u(t_k)=u'(\zeta_k)=0,\quad\forall\,k.$
\item[iii.] $u$ is increasing strictly on $(\zeta_{2k-1},\zeta_{2k})$ and decreasing strictly on $(\zeta_{2k},\zeta_{2k+1})$, $\quad\forall\,k$.
\end{description}
\end{lemma}
\ {\it Proof.} Suppose by contrast that the situation occurs. Using the functional energy $E(r)$, it is straightforward that $|u(\zeta_k)|\downarrow1$. Observe next that for $r$ large enough and $k\in\mathbb{N}$ unique such that $\zeta_{2k}\leq\,r<\zeta_{2k+1}$ or $\zeta_{2k+1}\leq\,r<\zeta_{2k+2}$, we have
$E(\zeta_{2k})\leq\,E(r)<E(\zeta_{2k+1})$ or $E(\zeta_{2k+1})\leq\,E(r)<E(\zeta_{2k+2})$ which means that $\displaystyle\lim_{r\rightarrow+\infty}E(r)=\displaystyle\frac{-\theta}{2(1-\theta)}$. In particular we get $\displaystyle\lim_{k\rightarrow+\infty}E(t_{k})=\displaystyle\frac{-\theta}{2(1-\theta)}$, which means that $\displaystyle\lim_{k\rightarrow+\infty}u'^2(t_{k})=\displaystyle\frac{-\theta}{1-\theta}<0$ which is a contradiction.

\end{document}